\input ./style/arxiv-general.cfg
\documentclass[aos,MSNbibl,nameyear,dvips]{arximspdf}
\makeatletter
   \@ifpackageloaded{graphicx}{}{\usepackage{graphicx}}
\makeatother
\usepackage{dcolumn}
\usepackage{graphicx}

\doi{10.1214/15-AOS1312}% Updated by VTEXPTS2LaTeX.exe, 27.02.2015
\volume{43}
\issue{3}
\pubyear{2015}
\firstpage{1323}
\lastpage{1350}
\docsubty{FLA}

\makeatletter
\newcolumntype{d}[1]{D{.}{.}{#1}}
\newcommand{\Prob}{\mathbb{P}}
\newcommand{\indep}{\perp\!\!\!\!\perp}

\newtheorem{thmm}{Theorem}[section]
\newtheorem{prop}{Proposition}[section]
\newtheorem{claim}{Claim}[section]
\newtheorem{lem}{Lemma}[section]
\newproclaim{rem}{Remark}[section]
\makeatother

\begin{document}
\begin{frontmatter}

\title{Higher criticism: $p$-values and criticism}
\runtitle{Higher criticism: $p$-values and criticism}

\begin{aug}
% Corresponding author: David Siegmund - siegmund@stat.Stanford.EDU% Updated by VTEXPTS2LaTeX.exe, 27.02.2015 16:06
%Updated by VTEXPTS2LaTeX.exe, 27.02.2015 10:20
\author[A]{\fnms{Jian} \snm{Li}\ead[label=e1]{jlijli@stanford.edu}}
\and
\author[A]{\fnms{David} \snm{Siegmund}\corref{}\ead[label=e2]{siegmund@stanford.edu}}
\runauthor{J. Li and D. Siegmund}
\affiliation{Stanford University}
\address[A]{Department of Statistics\\
Stanford University\\
Stanford, California 94305-4065\\
USA\\
\printead{e1}\\
\phantom{E-mail:\ }\printead*{e2}}
\end{aug}

% HISTORY:
%
\received{\smonth{1} \syear{2014}}% Updated by VTEXPTS2LaTeX.exe,
%27.02.2015 10:20
%
\revised{\smonth{10} \syear{2014}}% Updated by VTEXPTS2LaTeX.exe,
%27.02.2015 10:20

% ABSTRACT
%
\begin{abstract}
This paper compares the higher criticism statistic
(Donoho and Jin [\textit{Ann. Statist.} \textbf{32} (2004) 962--994]),
a modification of the higher criticism statistic also suggested
by Donoho and Jin, and two statistics of the Berk--Jones
[\textit{Z.~Wahrsch. Verw. Gebiete} \textbf{47} (1979) 47--59] type.
New approximations to the significance levels of the statistics are
derived, and their accuracy is studied by simulations.
By numerical examples it is shown that
over a broad range of sample sizes the Berk--Jones
statistics have a better power function than the higher
criticism statistics to detect sparse mixtures.
The applications suggested by
Meinshausen and Rice
[\textit{Ann. Statist.} \textbf{34} (2006) 373--393], to find
lower confidence bounds for
the number of false hypotheses, and by
Jeng, Cai and Li
[\textit{Biometrika} \textbf{100} (2013) 157--172], to detect copy number
variants, are also studied.
\end{abstract}

% KEYWORDS
% Pirmas kwd is didziosios raides
%
\end{frontmatter}

%s1 #&#
\section{Introduction}\label{sec1}
Donoho and Jin
(\citeyear{DJ})
consider the problem of deciding whether a large number, $n$, of
independently tested null hypotheses are all true, or whether some of them
are not
true. They discuss in detail
a suggestion of Tukey, called ``higher criticism,''
and they prove a number of asymptotic consistency results.
Suppose that $p_{(1)}, \ldots, p_{(n)}$ are ordered $p$-values for each
of the
individual hypotheses, which under the global null hypotheses that all
the individual null hypotheses are true, would be distributed as the
order statistics of a uniform sample on $[0,1]$. The test
statistics of Donoho and Jin are the higher criticism (HC) statistic
%
%e1 #&#
\begin{equation}
\label{hc} T_\mathrm{ HC} = n^{1/2} \max_{k_0 \leq k \leq k_1}
(k/n-p_{(k)})/ \bigl[p_{(k)}(1- p_{(k)})
\bigr]^{1/2},
\end{equation}
or a modified higher criticism statistic, which we will denote by
$T_{ \mathrm{ MHC}}$. This statistic is (\ref{hc}) modified by the
constraint that the $k$th term is included in
the statistic only if $p_{(k)} \geq1/n$.
The recommended values for the $k_i$ are $k_0 = 1$ and $k_1 = n/2$.
These statistics
reject the global null hypothesis if there is an excess of small
$p$-values.

Donoho and Jin find values $t_n$ such that
$\Prob_0 \{ T_\mathrm{ HC} \geq t_n \} \rightarrow0$, where $\Prob_0$ denotes
probability under the global null hypothesis, while under
certain ``borderline'' alternative configurations involving the
true number of nonnull hypotheses and a measure of their departures
from null, $\Prob\{ T_\mathrm{ HC} \geq t_n \} \rightarrow1$.
Thus the test of the global null is consistent, provided that there is
a certain minimal amount of separation between the global null hypothesis
and its negation. This is a pure hypothesis
testing problem, in the sense that the minimal amount of separation
is inadequate to allow one to identify with confidence which null
hypotheses are false, although one can be confident that they exist.

In studying this and related goodness-of-fit statistics based on
deviations of the empirical distribution function,
an approximate $p$-value based on
a classical result of Darling and
Erd\H{o}s (\citeyear{DE}) and adapted by \citet{J} is often cited.
However, since this approximation
is often very poor (see below), in practice $p$-values are often obtained by
simulation.

An alternative for relatively small sample sizes is the numerical
recursion of \citet{Noe},
which \citet{Owen} admirably exploited in finding confidence bands for a
distribution function.
\citet{ENS} give a substantially more efficient
algorithm, which requires $O(n^2)$ operations rather than the $O(n^3)$
required by No\'{e}'s.

The goals of this paper are: (i) to give approximations for the $p$-value of
the higher criticism type statistics that are reasonably accurate, even in
the situation of small $p$-values and large samples, where
numerical methods become onerous; (ii) to compare the power of a small
number of
different statistics that have the same basic properties of
consistency outlined above, but may perform differently in practice;
(iii) to illustrate application of our results by a discussion of two
papers that have developed related ideas for specific scientific
problems [viz. \citet{Rice}; \citet{LargeExp}].

The higher criticism statistic is suggested by standardizing $p_{(k)}$
as if it were asymptotically normally distributed even for small $k$.
As we shall see, this can exact a considerable price on the power of
the higher criticism statistic,
except when the number of false
null hypotheses is very small; this problem becomes very severe
when the significance level is small.

We have found particularly appealing an alternative class of statistics
suggested by \citet{BJ} as goodness of fit statistics, defined by
%
%e2 #&#
\begin{eqnarray}
\label{bj} T_\mathrm{ BJ} &=& \max_{k_0 \leq k \leq k_1} (2n)^{1/2}
\bigl\{(k/n) \log (k/np_{(k)})
\nonumber
\\[-8pt]
\\[-8pt]
\nonumber
&&\hspace*{70pt}{} + (1-k/n) \log \bigl[(1-k/n)/(1-
p_{(k)}) \bigr] \bigr\}^{1/2}.
\end{eqnarray}
For use in the context of higher criticism, we are interested in a one-sided
version of (\ref{bj}) where each term is modified by the condition that
$p_{(k)} < k/n$, which
we henceforth assume.
As explained below, a slightly
modified version designed to focus on
the small order
statistics is
%
%e3 #&#
\begin{eqnarray}
\label{mbj} T_\mathrm{ MBJ} &=& \max_{ k_0 \leq k \leq k_1} (2n)^{1/2}
I\{p_{(k)} < k/n\}
\nonumber
\\[-8pt]
\\[-8pt]
\nonumber
&&{}\times\bigl[(k/n) \log(k/np_{(k)}) - (k/n -
p_{(k)}) \bigr]^{1/2},
\end{eqnarray}
where the indicator function insures
focus on an excess of small $p$-values.

While \citet{BJ} suggest their statistic on the basis of
consideration of Bahadur efficiency, our preferred motivation is
in terms of Poisson variation. For example, suppose we observe a Poisson
process on $[0,1]$ and use the log likelihood ratio statistic to test
the hypothesis that the intensity is equal to one throughout the interval
against the alternative that there is a change-point at $t$, and
the intensity on $[0,t]$ is greater than one. Given that there
are $n$ events in the Poisson process, the
generalized likelihood ratio statistic
observed at the (ordered) times of the events is (\ref{bj}),
where now $p_{(k)}$ denotes the time of the $k$th event.
Alternatively, consider the empirical distribution
function $F_n(x)$ for a sample of size $n$ from the distribution
$F(x)$. For small $x$, $F_n(x)$ behaves for large $n$ like a
nonhomogeneous Poisson process, having log likelihood
$F_n(x) \log[F(x)] - F(x)$, which is maximized with respect to $F(x)$ at
$F_n(x) \log[F_n(x)] - F_n(x)$. To compare the lower tail of
the empirical distribution with the uniform
distribution, we consider the likelihood ratio
statistic $\max_x \{F_n(x) \log[F_n(x)/x] - [F_n(x) - x]\}$.
At the order statistics this becomes
(\ref{mbj}). When symmetrized by consideration of
both upper and lower tails, we get back to (\ref{bj}).

\citet{Walther} gives a similar argument in favor of the Berk--Jones
statistics.
(He then compares the Berk--Jones statistics with a completely
different class of statistics, the
``average likelihood ratio'' statistics, which seem to have
excellent power, but do not appear to be useful when estimation
of the number of
nonnull distributions [\citet{Rice}] or identification of
them is also desirable.)
In the context of goodness of fit, \citet{Jager} provide
asymptotic theory for a large class of
statistics, including those considered here, but they do not consider
the behavior of these
statistics for finite sample sizes.
Our methods apply to
many of these statistics, and in Section~\ref{sec5} we discuss briefly the one that
Jager and Wellner single out as perhaps a good
compromise between statistics that behave well in the tails
and those focusing on the center of the distribution.

The organization of the paper is as follows. In Section~\ref{sec2}, we give expressions
for approximate $p$-values, a heuristic argument in support of the
approximations and some Monte Carlo results demonstrating their accuracy.
In Section~\ref{sec3} we discuss comparative power.
In Section~\ref{sec4}, we revisit within our broader framework the research
of Meinshausen and Rice (\citeyear{Rice}), who discuss lower confidence bounds
on the number of false null hypotheses, also in a borderline case
where one cannot say exactly which null hypotheses are false.
We also re-examine briefly the interesting potential application
suggested by Jeng, Cai and Li (\citeyear{LargeExp}).
Proofs are given in the \hyperref[app]{Appendix}.

%s2 #&#
\section{Approximations}\label{sec2}

%s2.1 #&#
\subsection{Analytic approximations}\label{sec2.1}
We begin with the specific case of (\ref{hc}),
for which the calculations are more explicit, and then
suggest the minor modifications required for the other statistics.
The \hyperref[app]{Appendix} contains rigorous and detailed proofs for
the case of $T_{\mathrm{HC}}$ and $T_{\mathrm{MBJ}}$. The same argument with some technical
augmentation works for $T_{\mathrm{BJ}}$. The modified higher criticism
statistic requires still slightly different arguments, which
lead to slightly different approximations, given below and
in the \hyperref[app]{Appendix}.

Let $U_{(k)}$, $ k = 1, \ldots, n$ denote the order statistics
for a sample of size $n$
from a uniform $[0,1]$ distribution. For $1 \leq k_0 < k_1 < n$, let
\[
Z_n = \max_{k_0 \leq k \leq k_1} n^{1/2}[k/n -
U_{(k)}]/ \bigl[U_{(k)}(1-U_{(k)})
\bigr]^{1/2}.
\]
Let $C(x) = C(x,\xi) =
\{x+[\xi^2-\xi(\xi^2+4(1-x)x)^{1/2}]/2\}/(1+\xi^2)$, and
observe that
$Z_n \geq b$ if and only if $U_{(k)} \leq C(k/n, b/n^{1/2})$ for
some $k_0 \leq k \leq k_1$.
Hence the problem of approximating
$\Prob\{Z_n \geq b \}$
reduces to computing the sum over
$k_0 \leq k \leq k_1$ of the probabilities of the disjoint events
%
%e4 #&#
\begin{equation}
\label{basic1} \Prob \bigl\{U_{(k)} \leq C(k/n), U_{(k+j)} > C
\bigl((k+j)/n \bigr) \mbox{ for all } 1 \leq j \leq k_1 - k \bigr\}.
\end{equation}

The distribution of $U_{(k)}$ is
Beta with parameters
$k$ and $n-k+1$.

For our approximation we assume that $n \rightarrow\infty$.
From the joint distribution of
$U_{(i)}, i = k, \ldots, k_1$, it is easy to show by
calculation that the
joint conditional distribution of
$n[U_{(k+j)} - U_{(k)}]$ given that $U_{(k)} \sim C(k/n)$ converges to
the joint distribution of $\Gamma_j, j = 1, \ldots, k_1-k$, where
$\Gamma_j$ is the $j$th partial sum of independent, identically distributed
exponential random variables scale parameter $\lambda= (1-k/n)/(1-C(k/n))$.

Let $C'(x) = 1/(1+\xi^2)-\xi(1-2x)/\{(1+\xi^2)[\xi^2+4x(1-x)]^{1/2}
\}$
denote the derivative of $C(x) = C(x,\xi)$ with respect to $x$.
Given $U_{(k)} = C(k/n,b/n^{1/2}) - y/n$, the
conditional probability that $U_{(i)} > C(i/n,b/n^{1/2})$
for all $k < i \leq k_1$ converges to
%
%e5 #&#
\begin{eqnarray}
\label{rwprob} &&\Prob \bigl\{ \Gamma_j > j C'
\bigl(k/n,b/n^{1/2} \bigr) + y \mbox{ for all } 1 \leq j \leq
k_1-k \bigr\}
\nonumber
\\[-8pt]
\\[-8pt]
\nonumber
&&\qquad\sim \bigl[1-\lambda C'(k/n) \bigr]\exp(-
\lambda y),
\end{eqnarray}
provided $\lambda C'(k/n) < 1$, which will be the case if $C$ is
convex, as we assume throughout.

Let $c = C(k/n)$. Direct analysis of the probability
density function of $U_{(k)}$ shows that
\[
\Prob\{ U_{(k)} \in c - dy/n \} \sim f(c;k,n-k+1) \exp \bigl[-(k/n -
c)y/c(1-c) \bigr] \,dy/n,
\]
where $f$ denotes the Beta probability density function with the
indicated parameters.
Integrating asymptotically over
$[0,C(k/n,b/n^{1/2})]$
leads to the approximation for the term indexed by $k$,
%
%e6 #&#
\begin{equation}
\label{approx0} f(c; k,n+1-k) (c/k) \bigl[1-(1-k/n) c'/(1-c) \bigr],
\end{equation}
where $c = C(k/n,b/n^{1/2}), c' = C'(k/n,b/n^{1/2})$ and $f(x;\alpha,\beta)$
denotes the Beta density with parameters $\alpha,\beta$.
Our final approximation results from summing
(\ref{approx0}) over $k$.

Approximations for the other statistics involve obvious modifications.
For the Berk--Jones and modified Berk--Jones statistics, the curve
$C(x) = C(x,\xi)$ must be found numerically, while implicit
differentiation shows that $C'(x)$
is an explicit function of $C$ and $x$. For example,
for the modified Berk--Jones statistic, for which
$c = C(x) < x$ is the solution of $x \log(x/c) -(x-c) = \xi$,
by differentiation, we obtain
$C'(x) = \log[x/C(x)]/[1-x/C(x)]$. For the modified higher
criticism statistic we neglect terms where $1/n \geq C(k/n,b/n^{1/2})$,
and for other terms we modify the asymptotic value of the integral
over $[0,C(k/n,b/n^{1/2})]$ by subtracting from it
the asymptotic value of the integral over $[0,1/n]$.

\begin{rem*}The preceding argument is patterned after that of
\citet{Wood}, although the decomposition of a union of events
uses the last event that occurs rather than the first.
(Woodroofe's proof would be simplified by this change; here
it appears to be necessary.)
The \hyperref[app]{Appendix} gives a more detailed argument, which seems to be
unavoidably complex due to the fact that the events indexed by small
values of the subscript usually dominate the overall probability,
especially in the case of the higher criticism statistic.
Other techniques have been used to solve superficially similar problems,
but we were unable to use these. In particular,
we were unable to adapt the technique developed recently
by Yakir and colleagues to solve a variety of difficult
problems. See \citet{Yakir} and references cited there for many examples.
\end{rem*}
%s2.2 #&#
\subsection{Comparison with simulations}\label{sec2.2}

Table~\ref{table:DD} contains approximate $p$-values evaluated by
summing the terms of (\ref{approx0}) and comparison with simulations for
four different statistics: (i) the original higher criticism statistic (HC),
(ii) the modification
(MHC) obtained by requiring that $p_{(k)} \geq1/n$,
(iii) the (one-sided) Berk--Jones (BJ) statistic (\ref{bj}) and (iv)
the modification suggested in (\ref{mbj}) (MBJ).

In all cases $k_0 = 1$ and $k_1 = n/2$.
The number of repetitions of the Monte Carlo experiment is 100,000,
except in
the rows where $n = 30\mbox{,}000$, where it is 10,000.

Our approximations appear to be very good, although slightly
conservative, as one might conjecture from the derivations, which involve
approximating a convex curve by a sequence of successive tangents.

As might be anticipated, the significance
thresholds for the higher criticism statistic increase very rapidly
for decreasing significance levels. As we will
see below, this is the price that
the statistic pays to be able to detect very rare false null
hypotheses. But at very small significance levels, which are
appropriate for the application in Section \ref{sec4.2}, the threshold becomes
prohibitively large, unless one takes $k_0 > 1$, which calls into
question the advantage of HC for rare false null hypotheses.

%t1 #&#
\begin{table}
\caption{$p$-values}
\label{table:DD}
\begin{tabular*}{\textwidth}{@{\extracolsep{\fill}}ld{2.2}d{5.0}d{1.3}d{1.4}@{}}
\hline
\multicolumn{1}{@{}l}{\textbf{Statistic}} &
\multicolumn{1}{c}{\textbf{Threshold}} &
\multicolumn{1}{c}{$\bolds{n}$} &
\multicolumn{1}{c}{\textbf{Approximate} $\bolds{p}$\textbf{-value}} &
\multicolumn{1}{c@{}}{\textbf{Simulation}} \\
\hline
HC & 4.83 & 400 & 0.05 & 0.048 \\
HC & 10.0 & 400 & 0.01 & 0.01 \\
HC & 10.0 & 1000 & 0.01 & 0.010 \\
HC & 10.0 & 5000 & 0.01 & 0.010 \\
HC & 10.0 & 30\mbox{,}000 & 0.01 & 0.010 \\
HC & 31.0 & 1000 & 0.001 & 0.0009 \\
MHC & 3.13 & 400 & 0.05 & 0.053 \\
MHC & 3.91 & 400 & 0.01 & 0.010 \\
MHC & 3.94 & 1000 & 0.01 & 0.0101 \\
MHC & 3.98 & 5000 & 0.01 & 0.0098 \\
MHC & 4.00 & 30\mbox{,}000 & 0.01 &0.010 \\
MHC & 4.97 & 1000 & 0.001 & 0.0010 \\
BJ & 2.90 & 400 & 0.05 & 0.048 \\
BJ & 3.45 & 400 & 0.01 & 0.010 \\
BJ & 3.50 & 1000 & 0.01 & 0.0095 \\
BJ & 3.57 & 5000 & 0.01 & 0.0098 \\
BJ & 3.63 & 30\mbox{,}000 & 0.01 & 0.0096 \\
BJ & 4.14 & 1000 & 0.001 & 0.0009 \\
MBJ & 2.80 & 400 & 0.05 & 0.046 \\
MBJ & 3.35 & 400 & 0.01 & 0.0094 \\
MBJ & 3.40 & 1000 & 0.01 & 0.0094 \\
MBJ & 3.48 & 5000 & 0.01 & 0.0098 \\
MBJ & 3.56 & 30\mbox{,}000 & 0.01 & 0.0090 \\
MBJ & 4.04 & 1000 & 0.001 & 0.0009 \\
\hline
\end{tabular*}
\end{table}

%s2.3 #&#
\subsection{Other approximations}\label{sec2.3}
Diverse scientists writing about various aspects of the problems
considered in this paper and other related (often goodness of fit) problems
mention approximations based on
the double exponential extreme value distribution and attributed
to \citet{J}, who adapted the original result of
\citet{DE}. An intermediate step in deriving
this particular approximation involves the relation of the uniform empirical
process to a Brownian bridge, a step that makes the approximation
suspect, since the standard empirical distribution at small (or large) values
of its argument exhibits
Poisson, not Gaussian, variability.

Let $U(t)$ denote a stationary Ornstein--Uhlenbeck process with covariance
function $\exp(-|t|)$. Let $B_0(t)$ denote a Brownian bridge process
on $[0,1]$ and $W(t)$ a standard Brownian motion process on $[0, \infty)$.
Let $0 <\tau_0 < \tau_1 < 1$ and put $T_0 = 0.5 \log[\tau_1(1-\tau
_0)/\tau_0
(1-\tau_1)]$.
A Gaussian approximation that plays a role in the derivation of
the Darling--Erd\H{o}s result is given by
%
%e7 #&#
\begin{equation}\label{OUapprox}
\Prob \Bigl\{ \max_{\tau_0 \leq t \leq\tau_1} B_0(t)/ \bigl[t(1-t)
\bigr]^{1/2} \geq b \Bigr\} = \Prob \Bigl\{ \max_{0 \leq t \leq T_0 }
U(t) \geq b \Bigr\} \sim T_0 b \varphi(b),
\end{equation}
as $b \rightarrow\infty$, provided $T_0$ is bounded
or diverges slowly enough that the RHS of (\ref{OUapprox}) tends to 0. Here $\varphi$ is the standard normal
probability density function. For maxima over longer intervals, say
$[0,T]$, where
$T$ is so large that $T b \varphi(b)$ tends to a positive limit, which
we find it convenient to specify as $\exp(-x)$, it is easily shown by
consideration of $T/T_0$ approximately independent excursions of length $T_0$
that (\ref{OUapprox}) implies
$\Prob\{ \max_{0 \leq t \leq T } U(t) < b \} \rightarrow\exp(-\exp(-x))$.
If we put $T = \log(N)/2$, use the relationship between the
stationary Ornstein--Uhlenbeck
process and the Wiener process that\break
$U[\log(t)/2] = W(t)/t^{1/2}, 1 \leq t < \infty$, and
choose one particular asymptotic inversion of the limiting relationship
defined by
$ 0.5 \log(N) b \varphi(b) \rightarrow\exp(-x)$ to obtain $b$ as a
function of
$N$, we obtain the classical Darling--Erd\H{o}s approximation.
Most authors concede
that the approximation is very slow to converge.
\citet{Rice} suggest that the approximation only
be applied in the case of the modified higher criticism statistic,
and in this case it is not unreasonable. For example,
for the thresholds in rows 7--12 of Table~\ref{table:DD},
this approximation
would give 0.036,
0.008, 0.008, 0.008, 0.008 and 0.001.

For the same reasons that the Darling--Erd\H{o}s approximation seems to be
reasonable only in special cases, in particular for the modified higher
criticism statistic, the approximation (\ref{OUapprox}) is roughly correct
for the Berk--Jones statistic and slightly less so for the
modified Berk--Jones statistic, but not otherwise. For example, for
the thresholds in rows 13--18 of Table~\ref{table:DD}, this approximation
gives the $p$-values of 0.054, 0.01, 0.01, 0.01, 0.01, 0.001.
In fact, a different inversion of $b$ as a function of $N$ gives
the approximation suggested by \citet{WK} and
studied numerically by \citet{Walther}, which, as he shows, performs
reasonably well for
the Berk--Jones statistic, although not quite as well as direct
application of~(\ref{OUapprox}) in the cases we have tested.

It seems fair to say that there is enough latitude in performing this inversion
that one can frequently choose an approximation that seems to apply to
a particular problem.
It also seems clear from the results in Table~\ref{table:DD} that no
single approximation of this
type can be applied successfully to all four statistics, which while
asymptotically equivalent in the sense of \citet{DJ}, require
quite different thresholds to control their false positive probabilities.

%s3 #&#
\section{Power}\label{sec3}
In this section we consider the power under the commonly used
alternative model that the data arise from a mixture of a null
(often normal) distribution and a shifted version of the null
distribution.
%s3.1 #&#
\subsection{Analytic approximations}\label{sec3.1}
By a straightforward transformation
the evaluation of the power, that is, the probability under a mixture
model of rejecting the
global null, can be reduced to a problem of the same
structure as calculation of the global significance level, with one
important difference.
Suppose that $X_i$ are independent samples with
the distribution function $F(x)=(1-p)F_0(x) +
pF_0(x-\delta)$, where we consider in detail
the case $F_0 (x) = \Phi(x)$, the distribution function
of the standard normal distribution and $\delta> 0$. The
$p$-values are $p_j = 1-\Phi(X_j)$. We obtain independent
random samples from the uniform $[0,1]$ distribution by the following
transformation:
%
%e8 #&#
\begin{eqnarray}
\label{eqn:uni_tran} U_i & = & 1 - F(X_i)
\nonumber
\\
& = & 1 - F \bigl[\Phi^{-1}(1-p_i) \bigr]
\\
& = & (1-p) p_i + p \bigl\{ 1-\Phi \bigl[\Phi^{-1}(1-p_i)-
\delta \bigr] \bigr\}.\nonumber
\end{eqnarray}
Notice that $U_i$ is an increasing function of $p_i$, so the order
statistics of $U_i$ correspond directly to those of $p_i$ via equation
(\ref{eqn:uni_tran}). Let $d(i/n)$ denote
the transformed boundary $(1-p)C(i/n) + p \{
1-\Phi[\Phi^{-1}(1-C(i/n))-\delta]\}$. The global null is rejected
if and only if $U_{(k)}\leq d(k/n)$ for some $k_0\leq k\leq k_1$.

This new curve $d$ is not globally convex, so the argument of
Section~\ref{sec2.1} fails here. However, the curve is concave near 0,
and becomes convex after some point~$j_0$. Besides,
$\{U_{(i)}\}_{i<j_0}$ and $\{U_{(i)}\}_{i>j_0}$ are conditionally independent
given $U_{(j_0)}$. Therefore we have
\begin{eqnarray*}
& &\Prob \bigl\{U_{k}\leq d(k/n)\mbox{ for some }k_0
\leq k\leq k_1 \bigr\}
\\
&&\qquad =  \Prob \bigl\{U_{j_0}\leq d(j_0/n) \bigr\} + \int
_{d(j_0/n)}^1 f_{U_{(j_0)}}(x) \bigl[g_1(x)+g_2(x)-g_1(x)g_2(x)
\bigr]\,dx,
\end{eqnarray*}
where $g_1(x)=\Prob\{U_{(k)}\leq d(k/n)\mbox{ for some } k_0\leq
k<j_0|U_{(j_0)}=x\}$ and $g_2(x)=\Prob\{U_{(k)}\leq d(k/n)\mbox{ for
some } j_0<
k\leq k_1|U_{(j_0)}=x\}$. Then we approximate $g_2$ with the results in
Section~\ref{sec2.1}, and compute $g_1$ by No\'{e}'s recursion
[\citet{Noe}]. We observe fairly small $j_0$ for $(p,\delta)$ that gives
moderate values of power, leading to fast and
accurate implementation of  No\'{e}'s method, which encounters computational
difficulties for large $n$ if used by itself.

To illustrate the approximation above,
we consider the higher criticism and modified Berk--Jones
statistics for $n = 1000$, $\delta= 2.5$ and $p = 0.02$,
also for $\delta= 4$ and $p= 0.005$. The
significance level is 0.01. The values for
the power we obtained are
0.68 and 0.89 for higher criticism, and, respectively,
0.90 and 0.87 for the modified Berk--Jones statistic.
Simulations with 10,000 repetitions gave exactly
the same values to two significant figures.
It may be interesting
to note that the terms contributing substantially to the
power have indices considerably smaller than $np$, the
expected number of nonnull distributions. For the
higher criticism statistic, most of its power concentrates on the
first-order statistic, which is
the reason why it often performs poorly when the number of nonnull
distributions is not very small.

We have used  No\'{e}'s method because it is easy to understand and apply.
In a recent manuscript \citet{ENS}
describe an alternative, which after a number of
numerical refinements to improve its accuracy
appears to be
substantially faster [$O(n^2)$ instead of
$O(n^3)$ operations] and hence suitable for larger sample sizes.

%t2 #&#
\begin{table}
\caption{Power}
\label{table:PW}
\begin{tabular*}{\textwidth}{@{\extracolsep{\fill}}ld{5.0}d{2.2}d{1.1}d{1.3}c@{}}
\hline
\multicolumn{1}{@{}l}{\textbf{Statistic}} &
\multicolumn{1}{c}{$\bolds{n}$} & \multicolumn{1}{c}{\textbf{Threshold}} &
\multicolumn{1}{c}{$\bolds{\mu}$} & \multicolumn{1}{c}{$\bolds{p}$}& \multicolumn{1}{c@{}}{\textbf{Power}} \\
\hline
HC &400 & 4.83 &4.0 & 0.01 & 0.91 \\
BJ & & 2.90 & & & 0.87 \\
MHC & & 3.13 & & & 0.51 \\
MBJ & & 2.80 & & & 0.88 \\
HC &400 & 4.83 &1.5 & 0.1 & 0.54 \\
BJ & & 2.90 & & & 0.76 \\
MHC & & 3.13 & & & 0.73 \\
MBJ & & 2.80 & & & 0.79 \\
HC &1000 & 10.0 &1.5 & 0.08 & 0.19 \\
BJ & & 3.50 & & & 0.82 \\
MHC & & 3.94 & & & 0.81 \\
MBJ & & 3.40 & & & 0.81 \\
HC &1000 & 10.0 & 4.0 & 0.005 & 0.85 \\
BJ & & 3.50 & & & 0.81 \\
MHC & &3.94 & & & 0.43 \\
MBJ & & 3.40 & & & 0.83 \\
HC &1000 & 31.0 & 5.0 &0.002 & 0.67 \\
BJ & & 4.14 & & & 0.60 \\
MHC & & 4.97 & & & 0.04 \\
MBJ & & 4.04 & & & 0.62 \\
HC &1000 & 31.0 &2.0 & 0.05 & 0.11 \\
BJ & & 4.14 & & & 0.79 \\
MHC & & 4.97 & & & 0.78 \\
MBJ & & 4.04 & & & 0.80 \\
HC &5000 & 10.0 & 4.0 & 0.001 & 0.71 \\
BJ & & 3.57 & & & 0.65 \\
MHC & & 3.98 & & & 0.32 \\
MBJ & & 3.48 & & & 0.66 \\
HC &5000 & 10.0 & 3.0 & 0.003 & 0.53 \\
BJ & & 3.57 & & & 0.62 \\
MHC & & 3.98 & & & 0.55 \\
MBJ & & 3.48 & & & 0.63 \\
HC &5000 & 10.0 & 1.0 & 0.08 & 0.05 \\
BJ & & 3.57 & & & 0.81 \\
MHC & & 3.98 & & & 0.74 \\
MBJ & & 3.48 & & & 0.78 \\
HC &30\mbox{,}000& 10.0 & 3.0 & 0.001 & 0.51 \\
BJ & & 3.63 & & & 0.66 \\
MHC & & 4.00 & & & 0.65 \\
MBJ & & 3.56 & & & 0.68 \\
HC &30\mbox{,}000& 10.0 & 2.0 & 0.005 & 0.17 \\
BJ & & 3.63 & & & 0.68 \\
MHC & & 4.00 & & & 0.68 \\
MBJ & & 3.56 & & & 0.69\\
\hline
\end{tabular*}
\end{table}

%s3.2 #&#
\subsection{Power comparison}\label{sec3.2}
Now we compare by simulation the power of the four statistics discussed
above under the mixture model.
While the analytic method may offer computational advantages for very large
$n$, simulation has a number of compensating advantages for estimating
power, when the probabilities of interest are not small and a general
idea of their magnitude usually suffices, so very large sample
sizes are rarely required. One particular advantage in the case of
interest here is that statistical tests can be one-sided
or two-sided, and not only the number of nonnull distributions, but
also their noncentrality parameters can be variable
at essentially
no increase in computational effort.

In Table~\ref{table:PW} we compare by simulation the power of
the four statistics.
The listed thresholds correspond to examples in Table~\ref{table:DD}.
The alternative model is a mixture of $\mathrm{N}(0,1)$ and
$\mathrm{N}(\delta, 1)$ distributions.
The mixing parameter is $p$; the $\delta$'s are independent and have a
$\mathrm{N}(\mu, 0.1)$ distribution. The $p$-values are
two-sided. The number of repetitions of the simulation experiment
was 10,000, except for the last four scenarios, where it was reduced to 1000.
If instead of simulating a binomial($p$) number of nonnull distributions,
we take a deterministic $k = np$ number of nonnull distributions, the
power typically increases by roughly 10\%, except in some cases
when $p$ is very small. The general picture that emerges is that
at conventional levels of significance,
for very small $p$ the HC statistic can have about 5--8\% more power
than the two Berk--Jones statistics, which in turn have considerably
more power than the modified higher criticism statistic. For larger $p$,
the two Berk--Jones statistics and the modified higher criticism statistic
can have about equal power and substantially more power than the original
higher criticism statistic. At very small levels of significance, for
example, the level used in Section~\ref{sec4.2} where the higher criticism
statistic enters into a multiple comparisons analysis, the
HC statistic with the recommended $k_0 = 1$ can have very little power,
even for very small $p$.

%s4 #&#
\section{Applications}\label{sec4}
%s4.1 #&#
\subsection{Example: Confidence bounds for the proportion
of false null hypotheses}\label{sec4.1}

Consider a mixture model where $N\lambda$ of the hypotheses are false.
\citet{Rice}
give a lower confidence bound for $\lambda$,
which is based on a
functional of the process $[F_n(t)-t]/ \delta(t)$,
where $F_n(t)$ is the empirical distribution
function of the $p$-values, and $\delta(t)$ is a suitable function
chosen by the user. They
suggest the choice $\delta(t) = [t(1-t)]^{1/2}$, which
is closely related to the higher criticism statistic, and
they observe that the choice $\delta(t) = t$ is similarly related to
the Benjamini and Hochberg (\citeyear{BH}) false discovery rate criterion.
A similar lower confidence bound can be
obtained from the (modified) Berk--Jones statistic, which in view of
the power calculations of the preceding section, one might expect
to behave comparably or perhaps even better than the higher criticism statistic.
This subsection will
compare lower confidence bounds for the Gaussian mixture model.

We only consider the modified version, as
described in Section~\ref{sec1}, so in this
section we suppress the word modified in the description.

Assume that
(a) $n\gamma_{n,\alpha}$ is increasing in $n$, and that
(b) under the global null hypothesis $\Prob_0(\sup_{t\in\{s\in
(0,1)|F_n(s)\geq s
\}}[F_n(t)(\log F_n(t) - \log t)-(F_n(t)-t)]>\gamma_{n,\alpha})\leq
\alpha$ for
all $n$.

Define
\[
I_{\mathrm{BJ}} = \biggl\{\lambda\Big|\sup_{t:F_n(t)-\lambda\geq(1-\lambda)t}(F_n-
\lambda )\log\frac{F_n-\lambda}{(1-\lambda)t}- \bigl[F_n-\lambda-(1-\lambda )t
\bigr]> \gamma_{n,\alpha} \biggr\}
\]
and $\hat{\lambda}_{\mathrm{BJ}}=\sup I_{\mathrm{BJ}}$.
Then $\hat{\lambda}_{\mathrm{BJ}}$ can be shown to
be a lower confidence bound for $\lambda$ at
confidence level $1-\alpha$.

The proof of this result is similar to the argument given
by \citet{Rice} and hence is omitted.
To compute the required probability, we suggest using the
approximation obtained above. If $\gamma_{n,\alpha}$ is the
(approximate) $(1-\alpha)$
level quantile of the quantity in condition (b), the required
monotonicity condition (a) is satisfied numerically.

The simulation study reported below compares these two
lower confidence bounds. The underlying observations, $X_i,i=1,2,\ldots,N$,
are independently and normally distributed. A
fraction $\lambda N$ have a mean
of $\mu$ while
the others have a mean of 0. The lower confidence bounds, $\hat
{\lambda}_{\mathrm{MHC}}$ is calculated according to the
prescription of \citet{Rice}, and that for
$\hat{\lambda}_{\mathrm{MBJ}}$ is calculated according to the prescription
in the preceding paragraph.
The confidence level is
95\%. The ``bounding
sequence'' of Rice and Meinshausen, $\beta_{n,\alpha}$ and our
corresponding $\gamma_{n,\alpha}$ are
determined by the approximations given above.
Multiple configurations of model parameters, $N,\lambda$
and $\mu$, are considered,
and the simulations are repeated $10^5$ times for each configuration.

%t3 #&#
\begin{table}
\caption{Comparison between $\hat{\lambda}_{\mathrm{HC}}$
and $\hat{\lambda}_{\mathrm{BJ}}$, $N=400$}
\label{table:clbN400}
\begin{tabular*}{\textwidth}{@{\extracolsep{\fill}}lcd{1.4}d{1.2}cc@{}}
\hline
\multicolumn{1}{@{}l}{$\bolds{\lambda}$} & \multicolumn{1}{c}{$\bolds{\mu}$} &
\multicolumn{1}{c}{$\bolds{P(\hat{\lambda}_{\mathrm{HC}}>\hat{\lambda}_{\mathrm{BJ}})}$} &
\multicolumn{1}{c}{$\bolds{P(\hat{\lambda}_{\mathrm{HC}}<\hat{\lambda}_{\mathrm{BJ}})}$}
& \multicolumn{1}{c}{$\bolds{\|\lambda-\hat{\lambda}_{\mathrm{HC}}\|_2/\lambda}$}
&\multicolumn{1}{c@{}}{$\bolds{\|\lambda-\hat{\lambda}_{\mathrm{BJ}}\|_2/\lambda}$} \\
\hline
0.1 & 1.0 & 0.36 & 0.28 & 0.93 & 0.94 \\
0.1 & 1.5 & 0.43 & 0.55 & 0.79 & 0.89 \\
0.1 & 2.0 & 0.29 & 0.71 & 0.61 & 0.61 \\
0.1 & 2.5 & 0.26 & 0.74 & 0.46 & 0.45 \\
0.1 & 3.0 & 0.37 & 0.63 & 0.33 & 0.32 \\
0.2 & 1.0 & 0.59 & 0.40 & 0.80 & 0.83 \\
0.2 & 1.5 & 0.40 & 0.60 & 0.60 & 0.62 \\
0.2 & 2.0 & 0.17 & 0.83 & 0.44 & 0.44 \\
0.2 & 2.5 & 0.06 & 0.94 & 0.32 & 0.31 \\
0.2 & 3.0 & 0.03 & 0.97 & 0.23 & 0.21 \\
0.3 & 1.0 & 0.74 & 0.26 & 0.69 & 0.73 \\
0.3 & 1.5 & 0.45 & 0.55 & 0.50 & 0.51 \\
0.3 & 2.0 & 0.14 & 0.86 & 0.36 & 0.35 \\
0.3 & 2.5 & 0.03 & 0.97 & 0.26& 0.24 \\
0.3 & 3.0 & 0.001 & 0.99 & 0.18 & 0.16 \\
0.4 & 1.0 & 0.82 & 0.18 & 0.62 & 0.66 \\
0.4 & 1.5 & 0.42 & 0.58 & 0.43 & 0.44 \\
0.4 & 2.0 & 0.09 & 0.91 & 0.31 & 0.29 \\
0.4 & 2.5 & 0.02 & 0.99 & 0.22 & 0.19 \\
0.4 & 3.0 & 0.003 & 0.99 & 0.15 & 0.13 \\
0.5 & 1.0 & 0.84 & 0.16 & 0.57 & 0.61 \\
0.5 & 1.5 & 0.34 & 0.66 & 0.39 & 0.39 \\
0.5 & 2.0 & 0.04 & 0.96 & 0.27 & 0.25 \\
0.5 & 2.5 & 0.003 & 1.0 & 0.19 & 0.16 \\
0.5 & 3.0 & 0.0003 & 1.0 & 0.14 & 0.11 \\
\hline
\end{tabular*}
\end{table}

Numerical results for a sample size of $N=400$ are provided in Tables~\ref{table:clbN400} and \ref{table:clbl02}.
When the signal is weak, both methods can give a lower bound of 0, and
hence the sum of columns three and four can be less than one.
Two criteria are considered: (a) the larger of the
two lower confidence bounds,
and (b) the relative squared distances of the bounds from
the true parameter.
According to
the table, $\hat{\lambda}_{\mathrm{HC}}$ has an advantage over
$\hat{\lambda}_{\mathrm{BJ}}$ when the
values of $\mu$ are small. As $\mu$ increases, $\hat{\lambda}_{\mathrm{BJ}}$ first
becomes more precise in probability, and then it lies closer to
$\lambda$ than
$\hat{\lambda}_{\mathrm{HC}}$ does in (relative) $l_2$ distance. It should be
noted that even when $\Prob(\hat{\lambda}_{\mathrm{BJ}}>\hat{\lambda}_{\mathrm{HC}})$
is quite large, $\hat{\lambda}_{\mathrm{BJ}}$ may not be as
precise as $\hat{\lambda}_{\mathrm{HC}}$ in $l_2$ (e.g., $\lambda=0.1,
\mu=2$ and $\lambda=0.2, \mu=1.5$).
Therefore when $\lambda$ is small $\hat{\lambda}_{\mathrm{HC}}$ is better (in
probability and/or in $l_2$), while as $\lambda$ exceeds some critical
value $\lambda^*$, $\hat{\lambda}_{\mathrm{BJ}}$ becomes a tighter lower bound.
As can be seen in Table~\ref{table:clbN400}, the
borderline value of $\lambda^*$ for the probability comparison
seems to be rather stable for different values of
$\lambda$, whereas
the analogous $\lambda^*$ for $l_2$ distance decreases slowly from
above 2.0
to below 1.5 as $\lambda$ varies from 0.1 to 0.5. Moreover, when the
individual signal strength is weak ($\mu=1,1.5$), neither $\hat
{\lambda}_{\mathrm{HC}}$
nor $\hat{\lambda}_{\mathrm{BJ}}$ works well unless $\lambda$
is about 0.3 or even larger;
and in this case the difference between the two lower bounds
does not seem important compared to the
gap between the confidence bounds and the true $\lambda$.
This numerical behavior suggests
$\hat{\lambda}_{\mathrm{BJ}}$ is preferable unless either prior information
indicates a
weak individual signal in the data or the worst case scenario is
of primary concern.

%t4 #&#
\begin{table}
\caption{Comparison between
$\hat{\lambda}_{\mathrm{HC}}$ and $\hat{\lambda}_{\mathrm{BJ}}$,
$\lambda=0.2$}
\label{table:clbl02}
\begin{tabular*}{\textwidth}{@{\extracolsep{\fill}}lccccc@{}}
\hline
\multicolumn{1}{@{}l}{$\bolds{\lambda}$} & \multicolumn{1}{c}{$\bolds{\mu}$} &
\multicolumn{1}{c}{$\bolds{P(\hat{\lambda}_{\mathrm{HC}}>\hat{\lambda}_{\mathrm{BJ}})}$} &
\multicolumn{1}{c}{$\bolds{P(\hat{\lambda}_{\mathrm{HC}}<\hat{\lambda}_{\mathrm{BJ}})}$}
& \multicolumn{1}{c}{$\bolds{\|\lambda-\hat{\lambda}_{\mathrm{HC}}\|_2/\lambda}$}
&\multicolumn{1}{c@{}}{$\bolds{\|\lambda-\hat{\lambda}_{\mathrm{BJ}}\|_2/\lambda}$} \\
\hline
\phantom{0}400 & 1.0 & 0.59 & 0.40 & 0.80 & 0.83 \\
\phantom{0}400 & 1.5 & 0.40 & 0.60 & 0.60 & 0.62 \\
\phantom{0}400 & 2.0 & 0.17 & 0.83 & 0.44 & 0.44 \\
\phantom{0}400 & 2.5 & 0.06 & 0.94 & 0.32 & 0.31 \\
\phantom{0}400 & 3.0 & 0.03 & 0.97 & 0.23 & 0.21 \\
\phantom{0}800 & 1.0 & 0.70 & 0.30 & 0.71 & 0.75 \\
\phantom{0}800 & 1.5 & 0.49 & 0.51 & 0.51 & 0.53 \\
\phantom{0}800 & 2.0 & 0.20 & 0.80 & 0.37 & 0.37 \\
\phantom{0}800 & 2.5 & 0.06 & 0.94 & 0.27 & 0.25 \\
\phantom{0}800 & 3.0 & 0.02 & 0.98 & 0.19 & 0.17 \\
1200 & 1.0 & 0.81 & 0.19 & 0.66 & 0.71 \\
1200 & 1.5 & 0.60 & 0.40 & 0.47 & 0.49 \\
1200 & 2.0 & 0.25 & 0.75 & 0.33 & 0.33 \\
1200 & 2.5 & 0.07 & 0.93 & 0.24 & 0.22 \\
1200 & 3.0 & 0.03 & 0.97 & 0.17 & 0.15 \\
\hline
\end{tabular*}
\end{table}

%s4.2 #&#
\subsection{A more complex example}\label{sec4.2}
Motivated by the problem of detecting
intervals of copy number variation
(CNV) occurring at the same location in a (usually small)
fraction of aligned DNA sequences,
Jeng, Cai and Li (\citeyear{LargeExp}) suggest use of
a higher criticism based analysis as an alternative
to the method suggested by Zhang et al. (\citeyear{Zhang})
and \citet{Siegmund1}.
In brief, for each $n = 1, \ldots, N$, observations
$y_{n,t},t = 1, \ldots, T$ are independently and normally distributed
with constant known variances $\sigma_n^2$ and means that under
the null hypothesis are unknown constants $\mu_n$, but are different
by an increment $\delta_{n,I}$ in aligned short subintervals
$I \subset\{1,\ldots,T \}$. The subset of
$1, \ldots, N$ that exhibit changes in mean value in any particular
interval $I$
is usually relatively small. Jeng, Cai and Li's method
is, roughly speaking, to consider an interval
$I \subset\{1, \ldots, T \}$ having length at most $L$. They then
apply a higher criticism based analysis across the $N$ sequences to
a statistic (in this case the sample mean) defined on the interval $I$.
Large values of the higher criticism
on various intervals is interpreted as evidence that those intervals
contain CNV.

To control the false positive error rate,
they suggest using the approximation of Jaeschke referenced
above for each candidate interval in conjunction
with a Bonferroni bound (multiplication by $TL$) to account for
multiple comparisons involving overlapping candidate intervals
of different lengths.
For their actual analysis they use simulations.
The number of repetitions of their simulation experiments is
100 for a small set of data set and 50 for a larger set of data.

For this problem,
$N$ is often in the hundreds, $L$ is usually relatively small
while $T$ can be in the tens or hundreds of thousands.

Here we present a different simulation to compare a
higher criticism based procedure, along the lines
suggested by \citet{LargeExp}, a
modified higher criticism based procedure,
and its modified Berk--Jones counterpart.
The type I error is set to be approximately 0.05. The other parameters
are $N = 674$, $L = 20$
and $T = 40\mbox{,}929$.
Since the higher criticism statistic is extremely sensitive to the value
of $k_0$ in (\ref{hc}), we follow the suggestion
of \citet{LargeExp} and set $k_0 = 4$.

Although \citet{LargeExp} use this example to illustrate their
methods on real data,
for our comparative numerical experiment,
the data are similar in structure, but are artificially generated.
The number of intervals $I = [\tau_1, \tau_2]$ that contain signals
is 155,
of which 75
have a length of 3, 50 have a length of 4, 25 have a length
of 7 and
5 have a length of 10. The model has two variable parameters: given
that a particular interval contains at least one signal, $p$ is the
fraction of the $N$ intervals that contain the signal, and $\mu$ is
the change in mean values of the observed Gaussian random
variables in the interval and sequence that contain the signal.
The thresholds of the (modified) higher criticism
based procedure and the
(modified) Berk--Jones statistic are determined by
simulations repeated 900 times, and are compared
with our approximations. The number of repetitions in the power
computation is 625.

The significance thresholds obtained by simulation are
as follows (with theoretically calculated thresholds in parentheses):
(a) for a global false positive error rate of 0.05, HC $20.0(21.5)$,
MHC $9.3(9.1)$, MBJ $5.98(5.98)$;
(b) for a global false positive error rate
of 0.01, HC $24.1(26.0)$, MHC $9.84(9.79)$, MBJ $6.29(6.24)$.
Even though our theoretically
determined thresholds
are in principle conservative because of an inclusion of a Bonferroni
bound in the argument, the approximations in these examples appear to be
very good. This may not continue to be the case for larger values of
$L$.

Table~\ref{table:cmp3methods} shows the power of the three
procedures under different data configurations. Here power is taken
to be the
fraction of intervals containing signals that are detected.
Generally speaking the two higher criticism statistics have poor performance
for certain parameter combinations, small $p$ for MHC and not so small $p$
but small $\mu$ for HC, while the MBJ statistic maintains
good power througout the table.

%t5 #&#
\begin{table}
\caption{Power comparison}
\label{table:cmp3methods}
\begin{tabular*}{\textwidth}{@{\extracolsep{\fill}}lcccc@{}}
\hline
$\bolds{\mu}$ & $\bolds{p}$ & \textbf{Power using HC}
\textbf{(}$\bolds{k_0 = 4}$\textbf{)} & \textbf{Power using MHC} &
\textbf{Power using MBJ}\\
\hline
1 & 0.05& 0.11 & 0.16 & 0.17\\
1 & 0.06& 0.14 &0.21 & 0.22\\
1 & 0.07& 0.17 & 0.25 & 0.27\\
1 & 0.08& 0.19 &0.29 & 0.34\\
1 & 0.09& 0.21 &0.35 & 0.42\\
1.5 & 0.01& 0.12 &0.03 & 0.12\\
1.5 & 0.02& 0.26 & 0.19 & 0.28\\
1.5 & 0.03& 0.39 & 0.39 & 0.47\\
1.5 & 0.04& 0.53 & 0.60 & 0.67\\
2 & 0.01& 0.41 & 0.10 & 0.42\\
2 & 0.02& 0.79 & 0.55 & 0.81 \\
\hline
\end{tabular*}
\end{table}

As a final example, we compare the method of \citet{LargeExp},
based on the
modified Berk--Jones statistic, and
a method suggested without further study by Siegmund, Yakir
and Zhang (\citeyear{Siegmund1}), here denoted SYZ. Their method
contains a free parameter $p_0$, which
can be loosely interpreted as a prior expectation of the fraction of
the $N$ intervals that contain a signal whenever one is present.
Jeng, Cai and Li claim that their method is better at detecting
{\sl both} rare and common signals than a fixed value of $p_0$.
SYZ's
suggestion for making their method more robust against an incorrect
choice of $p_0$ was to use two different values, at say significance
level 0.025, so that by Bonferroni the overall significance level is $2
\times0.025 = 0.05$.
See also \citet{Xie}, who tested this suggestion in a somewhat
different context.
Here we consider the case where there is only a single interval $I$
containing a signal, which has various expected frequencies $p$ and
noncentrality parameters
$\xi= |I| \mu$ chosen so that the power is intermediate between 0 and 1.
The parameters $p_0$ of the SYZ procedure are chosen to equal
0.005 and 0.2, for which 0.025 significance thresholds are 25.3 and
179.2, respectively.
Table~\ref{table:cmpw/oldpaper} gives Monte Carlo estimates of the
marginal power, that is, the probability that
the statistic computed from observations from the correct interval $I$
exceed the appropriate
significance threshold.

%t6 #&#
\begin{table}
\caption{Power comparison}
\label{table:cmpw/oldpaper}
\begin{tabular*}{\textwidth}{@{\extracolsep{\fill}}ld{1.3}cc@{}}
\hline
\multicolumn{1}{@{}l}{$\bolds{\xi}$} & \multicolumn{1}{c}{$\bolds{p}$} &
\multicolumn{1}{c}{\textbf{Power using SYZ}} &
\multicolumn{1}{c@{}}{\textbf{Power using MBJ}}\\
\hline
5.0 & 0.005& 0.62 &0.52 \\
4.0 & 0.01 & 0.62 &0.47 \\
2.5 & 0.05 & 0.77 &0.53 \\
2.0 & 0.10 & 0.75 &0.56 \\
1.5 & 0.2 & 0.95 &0.75 \\
1.0 & 0.4 & 0.75 &0.56 \\
0.9 & 0.5 & 0.79 &0.62 \\
\hline
\end{tabular*}
\end{table}

%s5 #&#
\section{Discussion}\label{sec5}
We have derived an approximation to the significance level of higher
criticism like
statistics that appears to be sufficiently accurate for use in practice
and for
theoretical comparisons of the power of different statistics. As an
alternative to the
two higher criticism statistics suggested by \citet{DJ}, we
have also
studied two statistics motivated by the goodness-of-fit procedure
suggested by
\citet{BJ}. In a normal mixture
model, the Berk--Jones statistics have more power than the higher
criticism statistic,
except when the mixing fraction is very small, and more power than the modified
higher criticism statistic when the mixing fraction is small. Even in cases
where the Berk--Jones statistics have less power than one of the
higher criticism statistics, the
differences are only a few percent. The advantages of the Berk--Jones
statistics
are larger at smaller significance levels.
Since the significance threshold of the original higher criticism
statistic is extremely sensitive to the significance level,
when the test is an intermediate part of
a large multiple comparison problem (cf. Section~\ref{sec4.2}) and hence
involves a very small significance level,
its power can be much less than that of the other statistics.
This problem can be mitigated by taking a value $k_0 >1$ in definition
(\ref{hc}), but this deletes the capacity of the higher criticism statistic to
detect very rare mixtures. For the range of parameter values we have
studied, the two Berk--Jones statistics seem to be unequivocally better.

The statistics we have studied are related to goodness-of-fit tests
based on the empirical distribution function; but for the
higher criticism problem, as suggested by \citet{DJ} (and the
applications discussed in Section~\ref{sec4}),
we have focused on one-sided
statistics designed to detect an excess of small $p$-values.
\citet{Jager} develop an
elegant large
sample theory for a large class of statistics, but they do not
show how well their asymptotic theory predicts behavior
for sample sizes of interest in practice.

One statistic that receives particular mention by Jager and Wellner
as a perhaps reasonable compromise between statistics focusing on
the center of the distribution and statistics focusing on the tails
is (after modification to focus on an excess of small $p$-values)
$\max_{k_0 \leq j \leq k_1} n^{1/2}\{[(j/n)^{1/2} - p_{(j)}]^{+}\}^{1/2}$.
This statistic has the appealing feature that $C$ and $C'$ are given
explicitly by $C(x,\xi) = [(x^{1/2} - \xi)^{+}]^2$ and
$C'(x,\xi) = (1-\xi/x^{1/2})^{+}$.
Our methods apply and give good approximations
(compared to simulations) for the significance threshold.
For the examples in Table~\ref{table:PW}, we
find that the statistic behaves well for values of $p$ that are not too small.
It is usually more powerful than the modified higher criticism statistic,
but it has considerably less power than
the original higher criticism statistic and both Berk--Jones
statistics for small $p$. For example, for the third to fifth examples
in Table~\ref{table:PW}, we find by summation of (\ref{approx0}) that the threshold
$b = 1.54$ gives the same
level, 0.01 for $n = 1000$, as the examples given there, and we obtain
as estimates of the power
0.84, 0.57 and 0.27, respectively. For the
seventh to ninth examples, the appropriate threshold is 1.62, and
the power is 0.34, 0.47 and 0.82.

It might be interesting to see more systematically whether
our methods can be usefully applied in a goodness-of-fit context, for example,
as they might be applied to give confidence bands for a distribution function,
as in \citet{Owen}.

\begin{appendix}\label{app}
%s6 #&#
\section*{Appendix: Proofs}
The heuristic argument given above for our suggested approximations
is based on the approach of \citet{Wood} and uses
results obtained in a similar problem by \citet{Loader}.
Although the heuristic
is relatively simple, complete proofs are quite technical, and alternative
approaches that have been proved successful in apparently similar problems
do not seem to work here. The source of
the difficulties is the requirement that we not impose a lower bound on
$k_0$ and want $k_1$ to be of order $n$. In addition, different
statistics require somewhat different techniques. Here
we consider in detail the original higher criticism statistic
and the modified Berk--Jones
statistic.

We can obtain the analogous
approximations for the original Berk--Jones statistic by similar
methods (after
some additional technical arguments to verify the general conditions stated
below in Remark \ref{rem:conditions})
and for the Jager--Wellner statistic mentioned briefly in Section \ref{sec5}.
For the modified higher criticism statistic, we obtain by similar
calculations a slightly different
approximation given explicitly at the end of this appendix.

Consider the following two functions:
\begin{eqnarray*}
f_{1}(x,y) & = & \frac{x-y}{[y(1-y)]^{1/2}},
\\
f_{2}(x,y) & = & 2 \biggl[x\log\frac{x}{y}-(x-y) \biggr].
\end{eqnarray*}
For each $x,\xi$
let $C_i(x,\xi)$ denote the root $y \in(0,x)$ of $f_{i}(x,y) = \xi^i$.
Then
$C_{1}$ and $C_{2}$ correspond to the original higher criticism
statistic and modified Berk--Jones statistic, respectively.

Although $C_i$ is a function of two arguments, in most cases
the second argument will be $b/n^{1/2}$. When this is the case
we will simplify the notation by writing
$C_i(k/n)$.

%re6.1 #&#
\begin{rem}\label{rem:conditions}
We will see the proofs below hold in general for a piecewise
differentiable function $C(x), x\in[0,1]$ satisfying the following
conditions:
\begin{longlist}[(iii)]
\item[(i)] $C(x)$ is convex and $C(x)< x$ in the region of interest, $x=0$
excluded;

\item[(ii)] for some $\alpha\in(0,1)$, $C(k/n)\leq(k-1)/(n-1)$ for all $k
\in
[2,\alpha n]$ when $n$ is large enough;

\item[(iii)] for some $\alpha\in(0,1)$, $\sup_{0\leq x \leq
\alpha}(1-x)C'(x)/[1-C(x)] < 1$;

\item[(iv)] $\lim_{x\rightarrow0^{+}}C(x)/x=0$.
\end{longlist}
\end{rem}

The probabilities of rejecting the global null hypothesis with higher
criticism or Berk--Jones statistic have a similar expression.
%
%e9 #&#
\[
\Prob \Bigl(\max_{k_{0}\leq k\leq\beta
n}n^{1/2}f_1(k/n,p_{(k)})
\geq b \Bigr) =\Prob \Biggl(\bigcup_{k=1}^{\beta n}
\bigl\{p_{(k)}\leq C_{1}(k/n) \bigr\} \Biggr)
\]
and
%
%e10 #&#
\[
\Prob \Bigl(\max_{k_{0}\leq k\leq\beta n, p_{(k)}\leq
k/n}nf_2(k/n,p_{(k)})
\geq b^2 \Bigr) =\Prob \Biggl(\bigcup_{k=1}^{\beta n}
\bigl\{ p_{(k)}\leq C_{2}(k/n) \bigr\} \Biggr).
\]

If $k_0$ is proportional to $n$, the following division of the
rejection region is unnecessary, and Proposition~\ref{proposition:Loader} can be directly applied; otherwise the
rejection region should be divided into two
parts [equation (\ref{eqn:decomposition})], and their probabilities are
computed by different means.

There is an additional difficulty involving the value of $k_1$.
For our purposes $k_1 = n/2$ is the primary case of interest, and
in the following we take $k_1 = \beta n$ for a value $\beta< 1$.
In some cases, it is possible to take $k_1 = n-1$, and in others, this
imposes additional constraints on $k_0$; for instance, a sufficient
condition for higher criticism statistic is to have a constant $k_0$,
and in still others the
constraints are unclear.
%
%e11 #&#
\begin{eqnarray}
\label{eqn:decomposition}
& &\Prob \Biggl(\bigcup_{k=k_{0}}^{\beta n}
\bigl\{p_{(k)}\leq C(k/n) \bigr\} \Biggr)
\nonumber
\\
&&\qquad =  \Prob \Biggl(\bigcup_{k=k_0}^{\alpha
n-1} \bigl
\{p_{(k)} \leq C(k/n) \bigr\}{}\bigg\backslash{} \Biggl(\bigcup
_{k=\alpha
n}^{\beta n} \bigl\{p_{(k)}\leq C(k/n) \bigr\}
\Biggr) \Biggr)
\\
& &\qquad\quad{} +\Prob \Biggl(\bigcup_{k=\alpha n}^{\beta
n} \bigl
\{p_{(k)}\leq C(k/n) \bigr\} \Biggr).
\nonumber
\end{eqnarray}
The rejection regions of a large class of statistics, including the higher
criticism and Berk--Jones statistics, correspond to a collection of
curves $\{C(x)\}$, each of which satisfies $C(\beta)<\beta$ as well
as $C(0)=0$. Consequently there exists $\alpha\in(0,1/2)$ such that
$C(\beta)-C(\alpha)<\beta- \alpha$, and the $\alpha$ in
(\ref{eqn:decomposition}) satisfies these conditions.

Let $\Prob_{\mathrm{Bin}}(n,k,p)= {n \choose k }p^k(1-p)^{n-k}$ denote the
binomial probability distribution, $C'(x)=
\partial C(x,\xi)/ \partial x$. Proposition~\ref{proposition:Loader} below handles the second term of equation
(\ref{eqn:decomposition}).

%pr6.1 #&#
\begin{prop}\label{proposition:Loader}
Suppose that for every $\xi>0$, $C(x,\xi)$ is a convex and continuously
differentiable function of $x$, $C(x,\xi)$ is increasing in $x$ and
$C(x, \xi)<x$ for all $x\in[\alpha,\beta]$.
Then
\begin{eqnarray*}
& &\Prob \Biggl(\bigcup_{k=\alpha n}^{\beta n} \bigl
\{p_{(k)}\leq C(k/n) \bigr\} \Biggr)
\\
&&\qquad= \bigl(1+o(1) \bigr)\sum_{k=\alpha n}^{\beta n}
\biggl[1-\frac{(n-k+1)C'(k/n)}{n-nC(k/n)} \biggr] \Prob _{\mathrm{Bin}} \bigl(n,k,C(k/n) \bigr).
\end{eqnarray*}
\end{prop}

\begin{pf}
Let $F_{n}(x)$ be the empirical distribution function associated
with the independent $p$-values $p_{1},p_{2},\ldots, p_{n}$, and let
$D(x)$ be the inverse
of $C(x,\xi)$ with respect to $x$, that is, $D(C(x,\xi))=x$. Then
\begin{eqnarray*}
& & \Prob \Biggl(\bigcup_{k=\alpha n}^{\beta n} \bigl
\{p_{(k)}\leq C(k/n) \bigr\} \Biggr)
\\
&&\qquad = \Prob \bigl(F_{n}(x)\geq D(x)\mbox{ for some }x\in
\{p_{(\alpha
n)},p_{(\alpha n+1)},\ldots,p_{(\beta n)}\} \bigr)
\\
&&\qquad = \Prob \bigl(F_{n}(x)\geq D(x)\mbox{ for some }x\in \bigl[C(
\alpha, \xi),C(\beta,\xi) \bigr] \bigr)
\\
&&\qquad = \Prob \bigl(F_{n}(x)\leq1-D(1-x)\mbox{ for some }x\in \bigl[1-C(
\beta,\xi ),1-C(\alpha,\xi) \bigr] \bigr).
\end{eqnarray*}

The last equation results from the symmetry of $F_n(x)$, that is,\break $(\{
F_{n}(x)\}_{x\in[0,1]}\stackrel{d}{=}\{1-F_{n}(1-x)\}_{x\in[0,1]})$.
%Let $a(t)=1-D(1-t,\frac{b}{n^{1/2}})$, $\tau_{0}=1-C(\beta,b^*)$
%and $\tau_{1}=1-C(\alpha,b^*)$.
The desired result now follows from the
proof of Theorem~2.1 of \citet{Loader}.
\end{pf}

%re6.2 #&#
\begin{rem}
The summand in the formula in Theorem~2.1 of \citet{Loader} converges uniformly, so when
$\alpha n$ is not an integer, it could be replaced by $ \lceil
\alpha n \rceil$ or $ \lfloor\alpha n \rfloor$, and
the same goes for $\beta n$.
\end{rem}

The rest of the proofs show the first term on
the RHS of (\ref{eqn:decomposition}) has the identical
expression. The event in this term decomposes
into disjoint sub-events.
Let $B_{n,k} = \{p_{(k)}\leq C(k/n), p_{(k+j)}>C[(k+j)/n]
\ \forall j=1,2,\ldots,\beta n-k\}$.
Then this term equals $\sum_{k=1}^{\alpha n-1}\Prob(B_{n,k})$.

Let $f_{np_{(k)}}(x)$ denote the density of $np_{(k)}$, $f_{n,k}(y)$ be
$f_{np_{(k)}}(nC(k/n)-y)$
and $p_{n,k}(y)$ be $\Prob\{np_{(k+j)}>nC[(k+j)/n]\ \forall
j=1,2,\ldots,\beta n-k|np_{(k)}=nC(k/n)-y\}$.

%cl6.1 #&#
\begin{claim}
If $\{\varepsilon_{1},\varepsilon_{2},\ldots,\varepsilon_{n+1}\}$ is
a sequence of i.i.d. exponentially distributed random variables with mean
value of 1, and
$\Gamma_{k}=\sum_{i=1}^{k}\varepsilon_{i}$, then
$p_{n,k}(y)=\Prob(\Gamma_{j} / \Gamma_{n+1-k}>\frac
{y+nC[(k+j)/n]-nC(k/n)}{n+y-nC(k/n)}\ \forall j=1,\ldots,\beta n-k)$.
\end{claim}

\begin{pf}
The joint distribution of $(p_{(1)},\ldots,p_{(n)})$
is the same as that of $(\Gamma_{1}/\Gamma_{n+1},\ldots,\Gamma
_{n}/\Gamma_{n+1})$.
Conditional on $\Gamma_{k}/ \Gamma_{n+1}=nC(k/n)-y$,
$\{\Gamma_{k+j} /\break  \Gamma_{n+1} > nC[(k+j)/n]\}$ is
identical to $
\{(\Gamma_{k+j}-\Gamma_{k}) / (\Gamma_{n+1}-\Gamma_{k}) >\break
\frac{y+nC[(k+j)/n]-nC(k/n)}{n+y-nC(k/n)}\}$.

To complete the proof, we need to check the independence between
$\Gamma_{k} / \Gamma_{n+1}$
and
$\{(\Gamma_{k+j}-\Gamma_{k})/(\Gamma_{n+1}-\Gamma_{k})\}
_{j=1}^{n+1-k}$. Basu's
theorem indicates $(\Gamma_{n+1}-\Gamma_{k})\indep
\{(\Gamma_{k+j}-\Gamma_{k})/(\Gamma_{n+1}-\Gamma_{k})\}
_{j=1}^{n-k}$. Besides,
$\Gamma_{k}$ is independent of
$\{\Gamma_{k+j}-\Gamma_{k}\}_{j=1}^{n+1-k}$. Thus
$(\Gamma_{k},\Gamma_{n+1}-\Gamma_{k})\indep
\{(\Gamma_{k+j}-\Gamma_{k})/(\Gamma_{n+1}-\Gamma_{k})\}_{j=1}^{n-k}$,
which implies the desired independence.
\end{pf}

We know $p_{n,k}(y)$ is decreasing in $y$ for every
pair of $(n,k)$. The following claim shows $f_{n,k}(y)$ is
also decreasing in $y$ when $n$ is large enough.

%cl6.2 #&#
\begin{claim}\label{claim:monotonicity}
For all $n$ large enough and all $k=2,\ldots,\alpha n-1$, $f_{n,k}(y)$
is decreasing in $y$ when $C(x)$ is $C_{1}$ or $C_{2}$ [i.e., to check
condition \textup{(ii)} in Remark~\ref{rem:conditions}].
\end{claim}

\begin{pf}
Since $f_{np_{(k)}}$ is increasing on $[0,n(k-1)/(n-1)]$, it
suffices to show $C_{j}(k/n)<(k-1)/(n-1)$
for $j = 1,2$ when
$n$ is large enough. For any fixed $x$, $f_{1}(x,y)$ and $f_{2}(x,y)$
are decreasing in $y$ when $y\leq x$. Therefore the inequalities
are equivalent to $f_{j}(k/n,(k-1)/(n-1))<(b/n^{1/2})^j$
for $j = 1,2$, which results from the following limit, which converges
uniformly in $k$:
\begin{eqnarray*}
f_{1} \biggl(\frac{k}{n},\frac{k-1}{n-1} \biggr) & = &
\frac{1}{n}\sqrt{\frac{n-k}{k-1}}\rightarrow0\qquad (n\rightarrow\infty),
\\
f_{2} \biggl(\frac{k}{n},\frac{k-1}{n-1} \biggr) & = &
\frac{k}{n}\log\frac
{k(n-1)}{n(k-1)}-\frac{n-k}{n(n-1)}\rightarrow0\qquad (n
\rightarrow \infty).
\end{eqnarray*}
\upqed\end{pf}

%cl6.3 #&#
\begin{claim}\label{claim:onetermcuttail}
Assume that $b/n^{1/2}=\xi$ for some $\xi\in\mathbb{R}^{+}$, $C$ is
either $C_1$ or $C_2$ and let
$\delta=\log n$. Then there exists $ M=M(\xi)>1$ such that
\[
\frac
{\Prob(B_{n,k})}{\int_{0}^{nC(k/n)\wedge M\delta
}f_{n,k}(y)p_{n,k}(y)\,dy}\rightarrow1
\]
as $n\rightarrow\infty$ uniformly in $k$.
\end{claim}

\begin{pf}
Consider $g(x)=[x-C(x)]/\{C(x)[1-C(x)]\}$, $g$ is continuous on
$(0,\alpha]$; $\lim_{x\rightarrow0^+}g(x)=+\infty$ if
$\lim_{x\rightarrow0^+} C(x)/x=0$, a condition that $C_1$ and
$C_2$ satisfy. Hence $g$ achieves its minimum, denoted by $m$, on
$(0,\alpha]$. Let $M=1.1/m+1.1 >1$, and $I_{n,k}=\int_{0}^{nC(k/n)\wedge M\delta}f_{n,k}(y)p_{n,k}(y)\,dy$.
The claim reduces to
\begin{eqnarray*}
&&\int_{nC({k}/{n})\wedge
M\delta}^{nC({k}/{n})}f_{n,k}(y)p_{n,k}(y)\,dy\\
&&\qquad=
o(1)\int_{0}^{nC({k}/{n})\wedge
M\delta} f_{n,k}(y)
p_{n,k}(y) \,dy.
\end{eqnarray*}
Claim~\ref{claim:monotonicity} indicates that when $2\leq k\leq\alpha n-1$,
\begin{eqnarray*}
&&\frac{\Prob(B_{n,k})-I_{n,k}}{I_{n,k}}\\
&&\qquad \leq
\cases{ 0, &\quad  $\mbox{if $nC(k/n)\leq M
\delta$,}$ \vspace*{2pt}
\cr
\displaystyle\frac{[nC({k}/{n})-M\delta]^+
f_{n,k}(M\delta)p_{n,k}(M\delta)}{\delta f_{n,k}(\delta
)p_{n,k}(\delta)}, &\quad $\mbox{otherwise.}$}
\end{eqnarray*}

When $n\geq3$ and $nC(k/n)>M\delta$, we have $\delta=\log n>1$, and hence
%
%e12 #&#
\begin{eqnarray}
\label{equation:cuttail}
& & \frac{[nC(k/n)-M\delta]^{+}f_{n,k}(M\delta)p_{n,k}(M\delta
)}{\delta
f_{n,k}(\delta)p_{n,k}(\delta)}
\nonumber
\\[-8pt]
\\[-8pt]
\nonumber
&&\qquad = \frac{[nC(k/n)-M\delta]^{+}
[C(k/n)-M\delta/n]^{k-1}
[1-C(k/n)+M\delta/n]^{n-k}}{\delta[C(k/n)-\delta/n]^{k-1}
[1-C(k/n)+\delta/n]^{n-k}}.
\end{eqnarray}
Now consider the continuous version of the RHS of (\ref
{equation:cuttail}); that is, let
$k=nx, x\in(0,\alpha]$, $y=C(x)$. Recall that $\alpha$ is
less than $1/2$ in (\ref{eqn:decomposition}), so $y\leq x < 1/2$. When
$x$ satisfies $C(x)>M\delta/n>\delta/n$, we have
%
%e13 #&#
\begin{eqnarray}
\label{equation:continuous} & & \mbox{RHS of (\ref{equation:cuttail})}
\nonumber
\\
&&\qquad =  \frac{(ny-M\delta)^{+}}{\delta}\exp \biggl\{(nx-1)\log \biggl(1-\frac
{(M-1)\delta}{ny-\delta}
\biggr)
\nonumber
\\
& &\hspace*{80pt}\qquad\quad{}+n(1-x)\log \biggl(1+\frac{(M-1)\delta}{(1-y)n+\delta} \biggr) \biggr\}
\nonumber\\
&&\qquad\leq\frac{(ny-M\delta)^{+}}{\delta}\exp \biggl\{-(nx-1)\frac
{(M-1)\delta}{ny-\delta}+n(1-x)
\frac{(M-1)\delta}{(1-y)n+\delta} \biggr\}
\\
&&\qquad = \frac{(ny-M\delta)^{+}}{\delta}\exp \biggl\{-(M-1)\delta\frac
{(x-y)+{\delta}/{n}-{\delta}/{n^{2}}-
{(1-y)}/{n}}{(y-{\delta}/{n})(1-y+{\delta}/{n})} \biggr\}
\nonumber
\\
&&\qquad\leq\frac{(ny-M\delta)^{+}}{\delta}\exp \biggl\{-(M-1)\delta\frac
{(x-y)}{y(1-y)} \biggr\}\qquad
\mbox{for }n\mbox{ large enough.}\nonumber
\end{eqnarray}

The last inequality holds because $(y-\delta/n)(1-y+\delta/n)\leq y(1-y)$,
and $\delta/n-\delta/n^{2}-(1-y)/n$ is positive provided
$\log n>1+\log n/n$. Then the RHS of (\ref{equation:continuous})
$\leq\! ny/(\delta n^{1.1})\rightarrow0\ (n\rightarrow\infty).$
Hence the claim holds for all $k=2,3,\ldots,\alpha n-1$. For $k=1$ the
claim follows from $nC(1/n)\rightarrow0$.
\end{pf}

%le6.1 #&#
\begin{lem}\label{lemma:largedeviation}
If $\hat{\Gamma}_{n}=\sum_{i=1}^{n}(\xi_{i}-a)$
where $\xi_{i}$ are independent and
exponentially distributed with mean value of 1,
then $\limsup[\log\Prob(\hat{\Gamma}_{n}/n\in F)]/n=-\inf_{x\in
F}\Lambda_{a}^{*}(x)$
for any interval $F$ with positive length, where \textup{$\Lambda
_{a}^{*}(x)=a+x-1-\log(a+x)$
when $(a+x)>0$ and $=+\infty$ otherwise.}
\end{lem}

\begin{pf}
This lemma follows from the continuity of $\Lambda_{a}^{*}(x)$
and, for example, Theorem~2.2.3 in
\citet{Dembo},
page 27.
\end{pf}

In what follows we continue to use the notation introduced above:
$\delta=\log n$,
as in the condition of Claim \ref{claim:onetermcuttail},
$\Lambda_{a}^*$ is as described in Lemma \ref{lemma:largedeviation}
and $\Gamma_{k}=\sum_{i=1}^{k}\varepsilon_{i}$,
where the $\varepsilon_{i}$ are independent exponential random variables
with mean value~1.

%pr6.2 #&#
\begin{prop}\label{proposition:upper}
If the boundary function $C$ is $C_{1}$ or $C_{2}$, for all
$k=1,2,\ldots,\alpha n-1$,
then for any $\varepsilon>0$ we have $p_{n,k}(y) \leq(1+R_{n,k}(\varepsilon))
\{1-\frac{(n-k+1)C'(k/n)}{(1+\varepsilon) [n-nC(k/n)+M\delta]}\}
\exp\{-\frac{(n-k+1)y}{(1+\varepsilon) [n-nC(k/n)+M\delta]}\}$,
where $R_{n,k}(\varepsilon)\rightarrow0$ as $n\rightarrow\infty$
uniformly in $k$ and $y\in[0,nC(k/n)\wedge M\delta]$.
As a result\break $I_{n,k} \leq(1+ RR_{n,k}(\varepsilon))\int_{0}^{nC(k/n)\wedge
M\delta}f_{n,k}(y)
\{1-\frac{(n-k+1)C'(k/n)}{(1+\varepsilon)
[n-nC(k/n)+M\delta]}\}
\times\break \exp\{-\frac{(n-k+1)y}{(1+\varepsilon) [n-nC(k/n)+M\delta]}\}\,dy$
with $RR_{n,k}(\varepsilon)\rightarrow0$ uniformly in $k$.
\end{prop}

\begin{pf}
Since $C_{1}$ and $C_{2}$ are convex functions in $x$ for every $\xi$,
$C_{j}'$ is bounded above by 1 and bounded away from 1 when $x\leq
\alpha<
1/2$. Hence
%
%e14 #&#
\begin{eqnarray}
\label{equation:upperdetail} & & p_{n,k}(y)
\nonumber
\\
&&\qquad = \Prob \biggl(\frac{\Gamma_{j}}{\Gamma_{n+1-k}}\geq\frac
{y+nC((k+j)/n)-nC(k/n)}{n+y-nC(k/n)}\  \forall j
\in[1,\beta n-k] \biggr)
\nonumber
\\
&&\qquad \leq \Prob \biggl(\frac{\Gamma_{j}}{\Gamma_{n+1-k}}\geq\frac
{y+jC'(k/n)}{n+M\delta-nC(k/n)}\ \forall j
\in[1,\beta n-k] \biggr)
\nonumber
\\
&&\qquad \leq \Prob \biggl(\Gamma_{j}\geq\frac{(n+1-k)(y+jC'(k/n))}{(1+\varepsilon)
[n+M\delta-nC(k/n)]}\ \forall
j\in[1,\beta n-k] \biggr)
\\
& &\qquad\quad{} +\Prob \biggl(\biggl|
\frac{\Gamma_{n+1-k}}{n+1-k}\biggr|\in \bigl[(1+\varepsilon )^{-1},1+
\varepsilon \bigr]^{c} \biggr)
\nonumber\\
&&\qquad=  \Prob \biggl(\Gamma_{j}\geq\frac{(n+1-k)(y+jC'(k/n))}{(1+\varepsilon
)[n+M\delta-nC(k/n)]}\ \forall j
\geq1 \biggr)+\mathrm{Res}_{\mathrm{up}}
\nonumber
\\
& &\qquad\quad{} +\Prob \biggl(\biggl|
\frac{\Gamma_{n+1-k}}{n+1-k}\biggr|\in \bigl[(1+\varepsilon )^{-1},1+
\varepsilon \bigr]^{c} \biggr)
\nonumber
\\
&&\qquad =  \biggl\{1-\frac{(n-k+1)C'(k/n)}{(1+\varepsilon)[n-nC(k/n)+M\delta]}
\biggr\}\nonumber
\\
& &\qquad\quad{}\times\exp \biggl\{-\frac{(n-k+1)y}{(1+\varepsilon)[n-nC(k/n)+M\delta]} \biggr\} +\mathrm{Res}_{\mathrm{up}}
\nonumber
\\
& &\qquad\quad{} +\Prob \biggl(\biggl|\frac{\Gamma_{n+1-k}}{n+1-k}\biggr|\in \bigl[(1+\varepsilon )^{-1},1+
\varepsilon \bigr]^{c} \biggr).
\nonumber
\end{eqnarray}
The first term on the RHS of (\ref{equation:upperdetail}) is due to
8.13 of \citet{Sequential}, page 186.

Lemma~\ref{lemma:largedeviation} indicates
$\Prob(|\Gamma_{n+1-k}/(n+1-k)|\in[(1+\varepsilon)^{-1},1+\varepsilon
]^{c})\leq\break
A\exp\{-(1-\varepsilon)C_{\mathrm{up}}(\varepsilon)(n+1-k)\}$, where
$C_{\mathrm{up}}(\varepsilon)=\min\{\Lambda_{0}^{*}(1+\varepsilon),\Lambda
_{0}^{*}((1+\varepsilon)^{-1})\}>0$. The union bound of $\mathrm{Res}_{\mathrm
{up}}$ is
\begin{eqnarray*}
\mathrm{Res}_{\mathrm{up}} & \leq&\sum_{j=\beta
n-k}^{+\infty}
\Prob \biggl\{\Gamma_{j}<\frac
{(n+1-k)(y+jC'(k/n))}{(1+\varepsilon)
[n+M\delta-nC(k/n)]} \biggr\}
\\
& = & \sum_{j=\beta
n-k}^{+\infty}\Prob \biggl\{
\Gamma_{j}/j<\frac{[1-(k-1)/n]}{(1+\varepsilon)
[1+M\delta/n-C(k/n)]} \bigl[y/j+C'(k/n) \bigr]
\biggr\}.
\end{eqnarray*}
When $n$ is large enough $[1-(k-1)/n]C'(k/n)/\{(1+\varepsilon)
[1+M\delta/n-\break  C(k/n)]\}< a^*$ where
$a^{*}=\max_{x\in(0,\alpha]}(1-x)C'(x)/[1-C(x)]<1$. Since $j\geq
\beta n-k\geq(\beta-\alpha)n$ and $y\leq M\delta=
M\log n$, $[1-(k-1)/n]/\{(1+\varepsilon)
[1+M\delta/n-C(k/n)]\}y/j \leq\varepsilon' < 1-a^{*}$ for all $n$ large
enough. Let $C_{\mathrm{up}}^{*}(\varepsilon')$ denote
$\Lambda_{a^*}^*(\varepsilon')$, and hence Lemma~\ref
{lemma:largedeviation} provides the
upper bound of the summand,
\begin{eqnarray*}
\mathrm{Res}_{\mathrm{up}} & \leq& A\sum_{j=\beta n-k}^{+\infty}
\exp \bigl\{ -j(1-\varepsilon)C_{\mathrm{up}}^{*} \bigl(
\varepsilon' \bigr) \bigr\}
\\
& = &A\frac{\exp\{-(\beta-\alpha)n(1-\varepsilon)C_{\mathrm
{up}}^{*}(\varepsilon')\}}{1-\exp\{-(1-\varepsilon)C_{\mathrm
{up}}^{*}(\varepsilon')\}}.
\end{eqnarray*}
The second and third terms on the RHS of (\ref{equation:upperdetail})
decay uniformly faster than the first term, which tends to 0 more
slowly than $O(1/n)$.
\end{pf}

%pr6.3 #&#
\begin{prop}\label{proposition:lower}
With the same assumption of Proposition \ref{proposition:upper}, we have
$p_{n,k}(y)\geq(1+L_{n,k}(\varepsilon))[1-\frac{(1+\varepsilon
)(n-k+1)C'(k/n)}{n-nC(k/n)}]
\exp\{-\frac{(1+\varepsilon)(n-k+1)y}{n-nC(k/n)}\}$,\break 
where $L_{n,k}(\varepsilon)\rightarrow0$ as $n\rightarrow\infty$
uniformly in $k$ and $y\in[0,C_{n}(k/n)\wedge M\delta]$. As a result,
$I_{n,k} \geq
(1+LL_{n,k}(\varepsilon))\int_{0}^{nC(k/n)\wedge
M\delta}f_{n,k}(y) \{1-\frac{(1+\varepsilon
)(n-k+1)C'(k/n)}{n-nC(k/n)}\} \times
\exp\{-\frac{(1+\varepsilon)(n-k+1)y}{n-nC(k/n)}\}\,dy$
with $LL_{n,k}(\varepsilon)\rightarrow0$ uniformly in $k$.
\end{prop}

\begin{pf}
Due to the convexity of $C$, we have
%
%e15 #&#
\begin{eqnarray}
\label{eqn:lowerextra} & & p_{n,k}(y)
\nonumber
\\
&&\qquad=  \Prob \biggl(\frac{\Gamma_{j}}{\Gamma_{n+1-k}}\geq\frac
{y+nC((k+j)/n)-nC(k/n)}{n+y-nC(k/n)} \  \forall j
\in[1,\beta n-k] \biggr)
\nonumber
\\
&&\qquad\geq \Prob \biggl(\Gamma_{j}\geq\frac{(1+\varepsilon
)(n+1-k)}{n-nC(k/n)} \bigl[y+nC
\bigl((k+j)/n \bigr)-nC(k/n) \bigr]
\nonumber
\\
& & \hspace*{208pt}\qquad  \forall j \in[1,\beta n-k] \biggr) \nonumber\\
&&\quad\qquad{}-\Prob \biggl(\biggl|
\frac{\Gamma
_{n+1-k}}{n+1-k}\biggr| \in \bigl[(1+\varepsilon)^{-1},1+\varepsilon
\bigr]^{c} \biggr)
\nonumber
\\[-8pt]
\\[-8pt]
\nonumber
&&\qquad \geq \Prob \biggl(\Gamma_{j}\geq\frac{(1+\varepsilon)(n+1-k)}{n-nC(k/n)}
\bigl[y+jC' \bigl( \bigl(k+ \bigl\lceil n^{1/2} \bigr\rceil
\bigr)/{n} \bigr) \bigr]
\\
& &\qquad\hspace*{184pt}\forall j \in \bigl[1, \bigl\lceil n^{1/2} \bigr\rceil \bigr]
\biggr)\nonumber\\
&&\qquad\quad{}-\mathrm{Res}_{\mathrm
{down}}-\Prob \biggl(\biggl|\frac{\Gamma_{n+1-k}}{n+1-k}\biggr|\in
\bigl[(1+\varepsilon
)^{-1},1+\varepsilon \bigr]^{c} \biggr)
\nonumber
\\
&&\qquad\geq \Prob \biggl(\Gamma_{j}\geq\frac{(1+\varepsilon
)(n+1-k)}{n-nC(k/n)}
\bigl[y+jC' \bigl( \bigl(k+ \bigl\lceil n^{1/2} \bigr\rceil
\bigr)/{n} \bigr) \bigr]\ \forall j\geq1 \biggr)\nonumber
\\
& &\quad\qquad{} -\mathrm{Res}_{\mathrm{down}}-\Prob \biggl(\biggl|
\frac{\Gamma_{n+1-k}}{n+1-k}\biggr|\in \bigl[(1+
\varepsilon)^{-1},1+\varepsilon \bigr]^{c} \biggr).
\nonumber
\end{eqnarray}
By using the same argument as in Proposition \ref{proposition:upper}
and the uniform continuity of~$C'$,
\begin{eqnarray*}
& &\mbox{The first term on the RHS of (\ref
{eqn:lowerextra})}
\\
&&\qquad=  \biggl[1-\frac{(1+\varepsilon)(n+1-k)}{n-nC(k/n)}C' \bigl( \bigl(k+ \bigl\lceil
n^{1/2} \bigr\rceil \bigr)/n \bigr) \biggr]\\
&&\qquad\quad{}\times\exp \biggl\{-
\frac{(1+\varepsilon
)(n+1-k)}{n-nC(k/n)}y \biggr\}
\\
&&\qquad=  \bigl(1+o(1) \bigr) \biggl[1-\frac{(1+\varepsilon)(n+1-k)}{n-nC(k/n)}C'(k/n) \biggr]
\\
&&\qquad\quad{}\times\exp \biggl\{ -\frac{(1+\varepsilon)(n+1-k)}{n-nC(k/n)}y \biggr\}.
\end{eqnarray*}
Let $C'_{\max}=[C(\beta)-C(\alpha)]/(\beta-\alpha)<1$ [see
(\ref{eqn:decomposition})], we obtain $C((k+j)/n)-C(k/n)\leq jC'_{\max}/n$
for all $k<\alpha n$ and $k+j <\beta n$. Hence the union bound of
$\mathrm{Res}_{\mathrm{down}}$ becomes
\begin{eqnarray*}
& & \mathrm{Res}_{\mathrm{down}}
\nonumber
\\
&&\qquad\leq \sum_{j= \lceil n^{1/2} \rceil}^{\beta n - k}\Prob \biggl(
\Gamma_{j}\leq\frac{(1+\varepsilon
)(n+1-k)}{n-nC(k/n)} \bigl[y+nC \bigl((k+j)/n
\bigr)-nC(k/n) \bigr] \biggr)
\nonumber
\\
&&\qquad\leq \sum_{j= \lceil n^{1/2} \rceil}^{\infty}\Prob \biggl(
\Gamma_{j}\leq\frac{(1+\varepsilon)(n+1-k)}{n-nC(k/n)} \bigl(y+jC'_{\mathrm
{max}}
\bigr) \biggr)
\nonumber
\\
&&\qquad\leq A\sum_{j= \lceil n^{1/2} \rceil}^{\infty}\exp \bigl\{ -j(1-
\varepsilon)C_{\mathrm{down}}^{*}(\varepsilon) \bigr\}
\nonumber
\\
&&\qquad\leq A\frac{\exp\{- \lceil n^{1/2} \rceil(1-\varepsilon
)C_{\mathrm{down}}^{*}(\varepsilon)\}}{1-\exp\{-(1-\varepsilon)C_{\mathrm
{down}}^{*}(\varepsilon)\}},
\nonumber
\end{eqnarray*}
where $C_{\mathrm{down}}^{*}(\varepsilon)=\Lambda_{(1+\varepsilon
)a^{**}}^{*}(\varepsilon+\varepsilon^{2})>0$
with
$a^{**}=\max_{x\in(0,\alpha]}(1-x)C'_{\mathrm{max}}/[1 - C(x)]
$ (so $a^{**} < 1$).
By using an argument similar to that in Proposition \ref
{proposition:upper}, it can
be concluded that
$p_{n,k}(y)\geq(1+L_{n,k}(\varepsilon))[1-\frac{(1+\varepsilon
)(n-k+1)C'(k/n)}{n-nC(k/n)}]
\times\exp\{-\frac{(1+\varepsilon)(n-k+1)y}{n-nC(k/n)}\}$
with $L_{n,k}(\varepsilon)\rightarrow0$ uniformly in $k$.
\end{pf}

%cl6.4 #&#
\begin{claim}\label{claim:expfn}
Suppose that $f(y)$ is a nonnegative, nonincreasing function defined
on $[0,+\infty)$ with $f(a)>0$
for some $a>0$.
Then for any fixed $B>0$ and any $\beta_{1},\beta_{2}\geq B$, there
exists a continuous and increasing function, denoted by $\rho_B(x)$,
defined on $\mathbb{R}^+ \cup\{0\}$ with $\rho_B(0)=0$, such that
$|\ln\int_{\mathbb{R}^{+}}f(y)e^{-\beta_{1}y}\,dy-\ln\int_{\mathbb
{R}^{+}}f(y)e^{-\beta_{2}y}\,dy|\leq\rho_{B}(|\beta_{1}-\beta_{2}|)$
and that $\rho_B$ does not depend on $f$.
\end{claim}

\begin{pf}
Since
\begin{eqnarray*}
\int_{M_{1}}^{+\infty}f(y)e^{-\beta_{1}y}\,dy & \leq&
f(M_{1})e^{-\beta_{1}M_{1}}/\beta_{1}
\\
& = &e^{-\beta_{1}M_{1}}/ \bigl(1-e^{-\beta_{1}M_{1}} \bigr)\int_{0}^{M_{1}}f(M_{1})e^{-\beta_{1}y}\,dy
\\
&\leq& e^{-BM_{1}}/ \bigl(1-e^{-BM_{1}} \bigr)\int
_{0}^{M_{1}}f(y)e^{-\beta
_{1}y}\,dy,
\end{eqnarray*}
we have
\begin{eqnarray*}
\frac{\int_{\mathbb{R}^{+}}f(y)e^{-\beta_{1}y}\,dy}{\int_{\mathbb
{R}^{+}}f(y)e^{-\beta_{2}y}\,dy} & \leq&\frac{1}{1-e^{-BM_{1}}}\frac{\int_{0}^{M_{1}}f(y)e^{-\beta
_{1}y}\,dy}{\int_{0}^{M_{1}}f(y)e^{-\beta_{2}y}\,dy}
\leq\frac{e^{M_{1}|\beta_{1}-\beta
_{2}|}}{1-e^{-BM_{1}}},
\\
\bigl(M_{1}=|\beta_{1}-\beta_{2}|^{-{1}/{2}}
\bigr) & \leq&\frac
{e^{|\beta_{1}-\beta_{2}|^{0.5}}}{1-e^{-B|\beta_{1}-\beta
_{2}|^{-0.5}}}.
\end{eqnarray*}

Hence $\rho_{B}(x)=
|x|^{0.5}-\log(1-e^{-B|x|^{-0.5}})$ for $x\neq0$,
and = 0 otherwise is the desired function.
\end{pf}

Let $P_{n,k}$ denote the integral $\int_{0}^{nC(k/n)\wedge
M\delta}\{1-(n-k+1)C'(k/n)/[n-nC(k/n)]\}f_{n,k}(y)\exp\{-\frac
{n+1-k}{n-nC(k/n)}y\}\,dy$.

%pr6.4 #&#
\begin{prop}
As $n\rightarrow\infty$, $I_{n,k}/P_{n,k}\rightarrow1$ uniformly in $k$,
so consequently $\sum_{k=1}^{\alpha
n-1}\Prob(B_{n,k})\sim\sum_{k=1}^{\alpha
n-1}I_{n,k}\sim\sum_{k=1}^{\alpha n-1}P_{n,k}$.
\end{prop}

\begin{pf}
It is clear $1-\alpha\leq(n+1-k)/[n-nC(k/n)]\leq1$. Since
$C'(x)$ is bounded away from 1 when $0<x<\alpha$,
$1-\frac{(n-k+1)C'(k/n)}{(1+\varepsilon)[n-nC(k/n)+M\delta
]}=(1+O(\varepsilon)+O(M\delta/n))
\{1-\frac{(n-k+1)}{n-nC(k/n)}C'(k/n)\}
$.
For the same reason
$1-\frac{(1+\varepsilon)(n-k+1)}{n-nC(k/n)}\times C'(k/n)=(1+O(\varepsilon
))[1-\frac{(n-k+1)}{n-nC(k/n)}C'(k/n)]$.
So Claim \ref{claim:expfn} indicates
\begin{eqnarray*}
& &\int_{0}^{nC(k/n)\wedge M\delta}f_{n,k}(y)\exp \biggl
\{-\frac
{(n-k+1)y}{(1+\varepsilon)[n-nC(k/n)+M\delta]} \biggr\}\,dy
\\
&&\qquad\leq e^{\rho(O(\varepsilon)+O({M\delta}/{n}))}\int_{0}^{nC(k/n)\wedge M\delta}f_{n,k}(y)
\exp \biggl\{-\frac
{(n-k+1)y}{n-nC(k/n)} \biggr\}\,dy
\\
& &\qquad\quad{}\times \int_{0}^{nC(k/n)\wedge M\delta}f_{n,k}(y)\exp \biggl
\{-\frac
{(1+\varepsilon)(n-k+1)y}{n-nC(k/n)} \biggr\}\,dy
\\
&&\qquad\geq e^{\rho(O(\varepsilon))}\int_{0}^{nC(k/n)\wedge M\delta
}f_{n,k}(y)
\exp \biggl\{-\frac{(n-k+1)y}{n-nC(k/n)} \biggr\}\,dy.
\end{eqnarray*}
According to Propositions \ref{proposition:upper}, \ref
{proposition:lower}, the previous two inequalities
suggest $I_{n,k}/  P_{n,k}\rightarrow1 \mbox{ uniformly in
}k$ if we send $\varepsilon$ to 0 arbitrarily slowly.
\end{pf}

%pr6.5 #&#
\begin{prop}\label{proposition:firstterms}
The following convergence is uniform in $k$: $[1-\frac
{n+1-k}{n-nC(k/n)}C'(k/n)]\Prob_{\mathrm{Bin}}(n,k,C(k/n))/P_{n,k}\rightarrow
1$. As a result,
$\sum_{k=k_{0}}^{\alpha n-1}P_{n,k} \sim\sum_{k=k_{0}}^{\alpha
n-1}[1-\frac{n+1-k}{n-nC(k/n)}
C'(k/n)]\Prob_{\mathrm{Bin}}(n,k,C(k/n))$.
\end{prop}

\begin{pf}
(i) When $y\leq M\delta$ and $k\leq\alpha n$,
$\frac{y}{n-nC(k/n)}\leq M\delta/[n(1-\alpha)]\rightarrow0$,
so consequentially
\begin{eqnarray*}
& & \biggl[1-C(k/n)+\frac{y}{n} \biggr]^{n-k}\exp \biggl\{-
\frac{(n+1-k)y}{n-nC(k/n)} \biggr\}
\\
&&\qquad =  \bigl[1-C(k/n) \bigr]^{n-k}\exp \biggl\{(n-k)\log \biggl(1+
\frac{y}{n-nC(k/n)} \biggr)-\frac
{(n+1-k)y}{n-nC(k/n)} \biggr\}
\\
&&\qquad =  \bigl[1-C(k/n) \bigr]^{n-k}\exp \biggl\{-\frac{y}{n-nC(k/n)}-O
\biggl(\frac
{(n-k)y^{2}}{n^{2}[1-C(k/n)]^{2}} \biggr) \biggr\}.
\end{eqnarray*}
The remainder on the RHS of the last equation tends to 0 uniformly
in $k$ and $0\leq y\leq M\delta$. So
%
%e16 #&#
\begin{eqnarray}
\label{eqn:concentration1} P_{n,k} & \sim& \biggl[1-\frac{(n+1-k)C'(k/n)}
{n-nC(k/n)} \biggr]k\pmatrix{n \cr k} \bigl[1-C(k/n) \bigr]^{n-k}
\nonumber
\\[-8pt]
\\[-8pt]
\nonumber
& &{} \times\int_{0}^{nC(k/n)\wedge
M\delta} \biggl[C(k/n)-
\frac{y}{n} \biggr]^{k-1}\frac{1}{n}\,dy.
\end{eqnarray}
(ii)
The integral that appears on the RHS of (\ref{eqn:concentration1}) equals
%
%e17 #&#
\begin{equation}
\label{eqn:concentration2} \frac{1}{k} \biggl\{C \biggl(\frac{k}{n}
\biggr)^{k}- \biggl[C \biggl(\frac{k}{n} \biggr)-
\frac{nC(k/n)\wedge
M\delta}{n} \biggr]^{k} \biggr\},
\end{equation}
\begin{eqnarray*}
(\ref{eqn:concentration2})\Big/ \biggl[\frac{C(k/n)^{k}}{k} \biggr] & = &1- \biggl[1-
\frac{nC(k/n)\wedge
M\delta}{nC(k/n)} \biggr]^{k}
\nonumber
\\
& \geq& \cases{ 1, &\quad $\mbox{if }nC(k/n)<M\delta,$
\nonumber
\vspace*{2pt}
\cr
1-
\displaystyle\exp \biggl\{-\frac{kM\delta}{nC(k/n)} \biggr\}, &\quad
$\mbox{otherwise.}$}
\end{eqnarray*}
Since $k/[nC(k/n)]\geq1$,
$(\ref{eqn:concentration2}) = (1+o(1))
C(k/n)^{k}/k$ where the
infinitesimal tends to 0 uniformly in $k$.
\end{pf}

Propositions \ref{proposition:Loader}, \ref{proposition:firstterms}
together lead to the main result of this part.

%th6.1 #&#
\begin{thmm}[(Approximate formula for $p$-values of higher
criticism and Berk--Jones
statistics)] If the curve $C$ is $C_i$ ($i=1,2$) and $\beta\in(0,1)$,
then under the overall null
hypothesis
\begin{eqnarray*}
& &\Prob \Biggl(\bigcup_{k=k_0}^{\beta n} \bigl
\{p_{(k)}\leq C(k/n) \bigr\} \Biggr)
\nonumber
\\[-8pt]
\\[-8pt]
\nonumber
&&\qquad= \bigl(1+o(1) \bigr)\sum_{k=k_0}^{\beta n}
\biggl\{1-\frac{(n-k+1)C'(k/n)}{n[1-C(k/n)]} \biggr\} \Prob_{\mathrm{Bin}} \bigl(n,k,C(k/n)
\bigr).
\nonumber
\end{eqnarray*}
\end{thmm}

\begin{rem*} For the modified higher criticism statistic, the
approximation takes a slightly different form. The
binomial probability is replaced by
$\Prob_{\mathrm{Bin}}\{n,k, \max[1/n,C(k/n)]\} - \Prob_{\mathrm{Bin}}(n,k,1/n) \times
\max(nC(k/n),1)$.
\end{rem*}
\end{appendix}

\section*{Acknowledgments}
The authors would like to thank two referees and the Associate Editor
for their
careful reading and helpful suggestions to improve the exposition of
the paper.

% AOS,AOAS: If there are supplements please fill:
% \sname{Supplement A}
% \stitle{Title}
% \slink[doi]{10.1214/00-AOASXXXXSUPP}
% \sdatatype{.pdf}"
% \sdescription{Some text}

% imsref loaded by akundreckaite, 2015-02-27 10:42:14

%

% zodis "Acknowledgments" paliekamas pagal autoriu

%suskaldyti doi

\printaddresses
\end{document}